\documentclass[10pt]{article}

\usepackage{amsfonts,amsmath,amscd, amssymb,  vmargin}
\usepackage[francais]{babel}
\usepackage[all]{xy} 
\everymath{\displaystyle}
\begin{document} 
\title{Une neutralisation explicite de l'alg\`ebre de Weyl quantique compl\'et\'ee}
\author{Michel Gros\footnote{michel.gros@univ-rennes1.fr} \, et Bernard Le Stum\footnote{bernard.le-stum@univ-rennes1.fr}
\\IRMAR, UMR CNRS 6625\\Universit\'e de Rennes I\\Campus de Beaulieu\\35042 Rennes
cedex\\France }

\maketitle
\author

{\it{R\'esum\'e}}. Soit $p$ un nombre premier. Nous \'etablissons l'existence de  neutralisations  de  divers compl\'et\'es  de l'alg\`ebre de Weyl quantique sp\'ecialis\'ee en une racine de l'unit\'e primitive d'ordre $p$ (qui est ``g\'en\'eriquement'' une alg\`ebre d'Azumaya) et donnons en particulier un \'enonc\'e de neutralisation  explicite    relevant celui construit  en caract\'eristique $p$ dans \cite{gqlst}.

\setcounter{section}{-1}

\section{Introduction }

Si $X$ d\'esigne un sch\'ema lisse sur un corps parfait $k$ de caract\'eristique $p>0$, il est bien connu que le centre $Z{\cal{D}}_{X}^{(0)}$ de l'alg\`ebre des op\'erateurs de niveau $0$  (c'est-\`a-dire ceux ``sans puissances divis\'ees") ${\cal{D}}_{X}^{(0)}$ sur $X$ est ``gros'' et que ${\cal{D}}_{X}^{(0)}$ est une alg\`ebre d'Azumaya sur $Z{\cal{D}}_{X}^{(0)}$ (cf. par exemple \cite{OV} pour un rappel). Cette propri\'et\'e, raffin\'ee par la donn\'ee d'une neutralisation d'une compl\'etion centrale de ${\cal{D}}_{X}^{(0)}$,  donne lieu \`a diverses applications int\'eressantes (correspondance de Simpson en caract\'eristique positive \cite{OV}, classification d'\'equations diff\'erentielles \cite{vdp}, \'etude des repr\'esentations d'alg\`ebres de Lie modulaires \cite {bmr}....). Des variantes   de certaines de ces applications existent (e.g. \cite{ta}) dans un contexte ``quantique'' qui fournit donc un cadre possible pour rechercher   des rel\`evements modulo $p^{n}$ ($n\geq 1$) de la correspondance de Simpson en caract\'eristique $p$  d\'ecrite dans  \cite{gqlst}.  Dans cette note, nous donnons des neutralisations explicites de deux compl\'etions de l'alg\`ebre de Weyl quantique ${\text{D}}_{q}$   (sp\'ecialis\'ee en une racine primitive de l'unit\'e $q$ d'ordre $p$), qui est ``g\'en\'eriquement'' (cf. par exemple \cite{backelin}, Prop. 2.3) une alg\`ebre d'Azumaya. La plus int\'eressante d'entre elles   (cf. 4.1.3) rel\`eve celle obtenue en caract\'eristique $p>0$ dans \cite{gqlst}, Thm. 4.13 et est la clef de l'obtention des rel\`evements modulo $p^{n}$ auxquels on vient de faire allusion.  Dans un travail ult\'erieur, nous montrerons comment ce point de vue permet de reconsid\'erer certains d\'eveloppements  de la th\'eorie des \'equations aux $q$-diff\'erences.

\section{Notations}

Nous noterons $q$ une racine primitive de l'unit\'e  (dans une cl\^oture alg\'ebrique ${\overline{{\Bbb{Q}}}}$ de ${\Bbb{Q}}$ fix\'ee) d'ordre $p$ premier. Pour $n \geq 1$, on pose   $[n]=1+q+...+q^{n-1}$ et $[n]! = [n].[n-1]... [2].[1]$. Notons qu'on a   $[p]=0$.  

\section{D\'efinitions et Pr\'eliminaires}

On pose  $R := {\Bbb{Z}} [q]$ qu'on verra ci-dessous comme sous-anneau de ${\overline{{\Bbb{Q}}}}$. On notera $J$ l'id\'eal (maximal) de $R$ engendr\'e par $1-q$ et l'on identifiera $R/J$ \`a  $ {\Bbb{F}}_{p}$.  La classe  de $[n]$ modulo $J$ s'identifie alors simplement \`a la classe de $n$ modulo $p$ dans ${\Bbb{Z}}/p \simeq {\Bbb{F}}_{p}$.\\

{\bf{Lemme 2.1}}. (i) L'\'el\'ement $q$ est inversible dans $R$.

(ii) Les \'el\'ements $[i]$ avec  $(i,p)=1$ sont inversibles dans $R$. \\

{\it{D\'emonstration}}. (i) En effet, on a $q^{p}=1$.

(ii) Si $i$ est tel que $(i,p)=1$, on peut, dans ${\overline{{\Bbb{Q}}}}$,  \'ecrire    

(2.1.1)  \hspace{3cm} $\frac{1-q^{i}}{1-q}. \frac{1-q^{ij}}{1-q^{i}} = 1$ 
 d\`es que $ij \equiv 1 (p)$.

 Cette \'egalit\'e vaut en fait  dans $R$ et l'assertion en r\'esulte imm\'ediatement gr\^ace au th\'eor\`eme de B\'ezout.     $\square$

{\bf{D\'efinition 2.2}}.  On appelle {\it{alg\`ebre de Weyl quantique}} (sous-entendu sp\'ecialis\'ee en $q$)   le quotient ${\text{D}}_{q}$ de  l'anneau non commutatif $R \left\langle  x, \delta \right\rangle$ des polyn\^omes en deux variables $x$ et $\delta$    par  l'id\'eal engendr\'e par $\delta x- qx\delta-1$.

 L'\'el\'ement $\sigma := [\delta , x] := \delta x - x\delta =1-(1-q)x\delta \in {\text{D}}_{q}$ jouera un r\^ole important dans la suite.

  On peut faire agir $\delta \in {\text{D}}_{q}$ sur $R[x]$ en posant  $\delta(x^{n}) = [n]x^{n-1}$. Cette action s'\'etend en une action de ${\text{D}}_{q}$ sur $R[x]$ et l'on a $\sigma(x)=qx$ qui est donc un automorphisme de $R[x]$. On dispose ainsi d'un morphisme d'alg\`ebres ${\text{D}}_{q} \rightarrow {\text{End}}_{ {\Bbb{Z}}  }(R[x])$ (dont on prendra garde qu'il  n'est pas injectif : l'\'el\'ement $\delta^{p}$ a pour image l'endomorphisme nul). On notera que $\delta$ agit sur $R[x]$ comme une $\sigma$-d\'erivation,  au sens o\`u $\delta(fg) = \delta(f)g + \sigma(f)\delta(g)$ pour tous $f,g \in R[x]$. Plus bas interviendront diverses compl\'etions de ${\text{D}}_{q}$ et de $R[x]$ et nous noterons de la m\^eme mani\`ere les actions de $\delta$ et de $\sigma$.

 Soit ${\text{D}} = {\Bbb{F}}_{p}<x, \partial>/(\partial x- x \partial -1)$ l'alg\`ebre de Weyl (ordinaire) sur ${\Bbb{F}}_{p}$ qui n'est autre que l'alg\`ebre $\Gamma({\Bbb{A}}^{1}_{{\Bbb{F}}_{p}}, {\cal{D}}_{{\Bbb{A}}^{1}_{{\Bbb{F}}_{p}}}^{(0)})$ des sections globales du faisceau des op\'erateurs  diff\'erentiels de ``niveau $0$'' sur la droite affine ${\Bbb{A}}^{1}_{{\Bbb{F}}_{p}}$. 
  
{\bf{Proposition 2.3}}. L'application       

 (2.3.1) \hspace{3cm} ${\text{D}}_{q} \otimes_{R} {\Bbb{F}}_{p} \rightarrow {\text{D}}$ \,\,\, ;  \,\,\,$x\otimes 1 \mapsto x$,  $\delta \otimes 1 \mapsto  \partial$
 
 est un isomorphisme. \\
 
 {\it{D\'emonstration}}. En effet, le module ${\text{D}}_{q}$ (resp. ${\text{D}}$) est libre de base $\delta^{k}$ (resp. $\partial^{k}$) sur $R[x]$ (resp. ${\Bbb{F}}_{p}[x]$).$\square$

{\it{Remarque}}.  On pourrait (comme pour  l'alg\`ebre de Weyl ordinaire) introduire une th\'eorie d'alg\`ebre de Weyl quantique de niveau $m\geq0$ avec des puissances divis\'ees partielles de $\delta$ (cf. \cite{gqlst} pour le cas ordinaire) et d\'evelopper pour ces alg\`ebres une th\'eorie analogue \`a celle qui va suivre.

 \section{Sur la propri\'et\'e  Azumaya}

 Notons ${\text{Z}}_{q}$ le centre de l'alg\`ebre ${\text{D}}_{q}$ et $ {\text{Z}}(R[x])$ le centralisateur de $R[x]$ dans ${\text{D}}_{q}$. 
 
{\bf{Lemme 3.1}}. On a $[\delta, x^{p}]=[\delta^{p}, x]=[\delta^{p}, x^{p}]=0$. \\

{\it{D\'emonstration}}. C'est imm\'ediat \`a partir d'une des deux relations 

(3.1.1)   \hspace{3cm}      $\delta^{n}x=[n]\delta^{n-1} + q^{n}x\delta^{n}$ 
 
(3.1.2)      \hspace{3cm}     $\delta x^{n}=[n]x^{n-1}+q^{n}x^{n}\delta$  

obtenues par des r\'ecurrences sur $n$.\\

Soit $R[x^{p},\delta^{p}]$ (resp. $R[x,\delta^{p}]$) l'alg\`ebre commutative des polyn\^omes en les variables $x^{p}$ et $\delta^{p}$ (resp. $x$ et $\delta^{p}$).
 
{\bf{Proposition 3.2}}. (i) L'application canonique $R[x^{p},\delta^{p}] \rightarrow {\text{D}}_{q}$ induit  un isomorphisme 

(3.2.1)   \hspace{3cm}  $R[x^{p},\delta^{p}]  \stackrel \simeq  \longrightarrow  {\text{Z}}_{q}$.

(ii) L'application canonique $R[x,\delta^{p}] \rightarrow  {\text{D}}_{q}$ induit un isomorphisme 

(3.2.2)   \hspace{3cm} $R[x,\delta^{p}]  \stackrel \simeq  \longrightarrow  {\text{Z}}(R[x])$. \\
 
{\it{D\'emonstration}}. (i) Tout \'el\'ement $P \in {\text{D}}_{q}$ s'\'ecrit de mani\`ere unique $P= \sum_{i=0}^{n} a_{i}(x)\delta^{i}$ avec $a_{i}(x) \in R[x]$. Compte tenu de la relation $(3.1.1)$, la relation de commutation  $Px=xP$ implique, en examinant l'\'egalit\'e obtenue degr\'e par degr\'e en $i$ dans l'\'ecriture ci-dessus, que $i$ doit \^etre un multiple de $p$. Un argument analogue avec la commutation \`a $\delta$ (utilisant cette fois   $(3.1.2)$) donne l'appartenance de chaque $a_{i}(x)$ \`a $R[x^{p}]$. 

(ii) On proc\`ede de mani\`ere analogue. $\square$

Dans la suite, pour tout $R$-module $M$, on notera $\widetilde M$ son s\'epar\'e-compl\'et\'e $J$-adique (ou $p$-adique : c'est la m\^eme chose car $(p) \subset J$ et $J^{p} \subset (p)$).

{\bf{Proposition 3.3}}.  L'alg\`ebre $\widetilde{{\text{D}}_{q}}$  est une $\widetilde { {\text{Z}}_{q}}$-alg\`ebre d'Azumaya  de rang $p^{2}$ neutralis\'ee par l'extension $ \widetilde {{\text{Z}}_{q}} \rightarrow \widetilde {{\text{Z}}(R[x])}$.  \\

 {\it{D\'emonstration}}. Montrons en effet  que  l'application canonique    
 
 (3.3.1) \hspace{1cm}   $\widetilde{{\text{D}}_{q}} \otimes_{\widetilde {{\text{Z}}_{q}}} \widetilde {{\text{Z}}(R[x])} \rightarrow {\text{End}}_{\widetilde {{\text{Z}}(R[x])}}(\widetilde{{\text{D}}_{q}})$ \,\,\,; \,\,\, $Q \otimes Q' \mapsto (P \mapsto QPQ')$  
 
est un isomorphisme. Les modules en pr\'esence  \'etant de m\^eme rang $p^{2}$, il suffit d'\'etablir  la surjectivit\'e de (3.3.1). De plus, comme ils sont $J$-adiquement complets (car  de type fini sur  $\widetilde {{\text{Z}}(R[x])}$), il suffit donc, gr\^ace au lemme de Nakayama, de v\'erifier la surjectivit\'e de la r\'eduction modulo $J$ de (3.3.1).  Comme cette derni\`ere s'identifie   (avec les notations \'evidentes calqu\'ees des pr\'ec\'edentes) \`a 

(3.3.2) \hspace{1cm} $ {\text{D}}    \otimes_{{\text{Z}}}   { {\text{Z}}({\Bbb{F}}_{p}[x])} \rightarrow {\text{End}}_{  {\text{Z}}({\Bbb{F}}_{p}[x]) }({\text{D}})$ \,\,\,; \,\,\,  $Q \otimes Q' \mapsto (P \mapsto QPQ')$

la surjectivit\'e d\'ecoule    imm\'ediatement du r\'esultat de M. Kaneda rappel\'e dans \cite{gqlst}, thm. 3.7. $\square$

{\it{Remarques}}. (i) Rappelons (\cite{backelin}, lemma 2.4) qu'on a  $\sigma^{p} = 1-(1-q)^{p}x^{p}\delta^{p} \in {\text{Z}}_{q}$. On peut alors reformuler (une fois \'etendu ici les scalaires de $R$ \`a ${\overline{{\Bbb{Q}}}}$) \cite{backelin}, Prop. 2.3 en disant que ${\text{D}}_{q}$ est une  ${\text{Z}}_{q}$-alg\`ebre d'Azumaya au-dessus du lieu d'inversibilit\'e ${\cal{U}} \subset {\text{Spec}}\,({\text{Z}}_{q})$ de  $\sigma^{p}$.  
Le lecteur observera que passer \`a la compl\'etion $J$-adique rend automatiquement  $\sigma^{p}$ inversible (car $(1-q)=J$).

(ii) Nous ignorons si l'analogue \'evident de (3.3.1) avant $J$-compl\'etion est un isomorphisme au-dessus de ${\cal{U}}$ (ce qui raffinerait \cite{backelin}, Prop. 2.3). C'est, comme on va le voir tout de suite,  par exemple le cas  pour $p=2$  ($q=-1$), exemple qui va aussi nous servir \`a montrer o\`u intervient  l'inversibilit\'e de $\sigma^{p}$ dans la th\'eorie. En effet, si $\sigma^{2}$ est inversible, alors, dans la base $(1,\delta)$ de ${\text{D}}_{q}$ sur ${\text{Z}}(R[x])$, les quatre matrices  

(3.3.3) \hspace{1cm} $E_{1}:=\left(\begin{array}{cc}0 & \sigma^{2}\\0 & 0\end{array}\right)$, $E_{2}:=\left(\begin{array}{cc}0 & 0\\0 & \sigma^{2}\end{array}\right)$,  $E_{3}:=\left(\begin{array}{cc}\sigma^{2} & 0\\0& 0\end{array}\right)$, $E_{4}:=\left(\begin{array}{cc}0 & 0\\\sigma^{2}& 0\end{array}\right)$   

 forment une base   de ${\text{End}}_{{\text{Z}}(R[x])}({\text{D}}_{q})$  (au-dessus de ${\cal{U}}$) et l'isomorphisme cherch\'e s'obtient en v\'erifiant que

(3.3.4)   \hspace{1,9cm} $u:=(1\otimes x   - x\otimes 1) + 2 [x \otimes x\delta - x^{2} \otimes \delta]    \mapsto E_{1}$, 
  
   \hspace{3cm} $v:=  u. \delta \mapsto E_{2}$, 

   \hspace{3cm} $1\otimes \sigma^{2}-  v   \mapsto E_{3}$, 
   
   \hspace{3cm} $1\otimes \sigma^{2}\delta -1\otimes u \delta^{2} \mapsto E_{4}$.

\section{Neutralisation}

Soit  $I$ l'id\'eal  de ${\text{Z}}_{q}$ engendr\'e par $\delta^{p}$. Dans la suite, pour tout ${\text{Z}}_{q}$-module $M$, on notera $\widehat{M}$ son s\'epar\'e compl\'et\'e $I$-adique. Comme il r\'esultera de la proposition 4.1 ci-dessous,  le compl\'et\'e $\widehat{ {\text{Z}} _{q}}$ s'identifie bien au   centre de $\widehat{ {\text{D}}_{q}}$, ne cr\'eant ainsi aucune ambiguit\'e dans les notations. Par commodit\'e, nous utiliserons librement dans la suite la base $(1,x,\, .... \,,x^{p-1})$ du  $\widehat{ {\text{Z}} _{q}}$-module  $\widehat{ {\text{Z}}(R[x])}$. 
 
{\bf{Proposition 4.1}}. Soient $D = (d_{i,j})_{1\leq i,j \leq p}, X= (x_{i,j})_{1\leq i,j \leq p} \in  {\text{End}}_{\widehat{ {\text{Z}} _{q}}}(\widehat{ {\text{Z}}(R[x])})$   les endomorphismes d\'efinis (dans la base ci-dessus) par les matrices dont les seuls coefficients non nuls sont les suivants :

(4.1.1) \hspace{3cm} $d_{i,i+1}= [i]+ q^{i}x^{p}\frac{\delta^{p}}{[p-1]!}$ pour $i \in [1,p-1]$  et $d_{p,1} = \frac{\delta^{p}}{[p-1]!}$  ;

(4.1.2) \hspace{3cm} $x_{1,p}=x^{p}$ et  $x_{i,i-1} =1$ pour $i\in [2,p]$.  

Alors, l'application de $R$-alg\`ebres $R \left\langle  x, \delta \right\rangle \rightarrow {\text{End}}_{\widehat{ {\text{Z}}_{q}}}(\widehat{ {\text{Z}}(R[x])})$ d\'efinie par $\delta \mapsto D$ ; $ x  \mapsto X$ induit, par passage au quotient,   un isomorphisme de $\widehat{ {\text{Z}} _{q}}$-alg\`ebres

(4.1.3) \hspace{3cm}  $\widehat{ {\text{D}}_{q}}  \stackrel \simeq  \longrightarrow {\text{End}}_{\widehat{ {\text{Z}} _{q}}}(\widehat{ {\text{Z}}(R[x])}).$  

De plus, la r\'eduction modulo $p$ de (4.1.3) s'identifie canoniquement \`a l'isomorphisme du thm. 4.13 de \cite{gqlst}. \\

 {\it{D\'emonstration}}. On a donc  
 
$$  D=\left(\begin{array}{cccccc}
0 & [1]+qx^{p} \frac{\delta^{p}}{[p-1]!}  & 0 &  & \cdots & 0 \\ \\
\vdots & \ddots & [2]+q^{2}x^{p} \frac{\delta^{p}}{[p-1]!} && \ddots  & \vdots \\ \\
\vdots &  &  && \ddots & 0 \\ \\
0 &  & \ddots &&  & [p-1]+q^{p-1}x^{p} \frac{\delta^{p}}{[p-1]!}  \\ \\
\frac{\delta^{p}}{[p-1]!} & 0 & \cdots && 0 & 0
\end{array}\right)$$  \\ 
 
$$ X=\left(\begin{array}{ccccc}
 0 & \cdots  & \cdots & 0 & x^{p} \\
 1 & \ddots & &  & 0 \\
 0 & 1 & \ddots &  & \vdots \\
\vdots & \ddots & \ddots & \ddots & \vdots  \\
 0 & \cdots & 0 & 1 & 0\end{array}\right) \,\,\,.$$   
 
Les seuls coefficients non nuls de $DX$ (resp. $XD$) sont  situ\'es sur la diagonale et le $(j, j)$-i\`eme vaut  

(4.1.4) \hspace{3cm} $[j]+x^{p}q^{j}\frac{\delta^{p}}{[p-1]!}$ (resp.  $[j-1]+x^{p}q^{j-1}\frac{\delta^{p}}{[p-1]!}$ )  

pour tout $j\in [1,p]$ (avec ici l'abus d'\'ecriture $[0]=0$). La matrice $DX-qXD-1$ a donc pour $(j, j)$-i\`eme coefficient

(4.1.5) \hspace{3cm} $[j]+x^{p}q^{j}\frac{\delta^{p}}{[p-1]!} - q([j-1]+x^{p}q^{j-1}\frac{\delta^{p}}{[p-1]!} -1)$

qui est effectivement nul car $[j]=1+q[j-1]$ pour tout $j\in [1,p]$. 
 
 On dispose donc bien de l'application (4.1.3). On va ensuite utiliser un argument de r\'eduction modulo $I$.  Etablissons tout d'abord le 
 
 {\bf{Lemme 4.2}}. L'action naturelle de ${\text{D}}_{q}$ sur $R[x]$ est $R[x^{p}]$-lin\'eaire et induit un isomorphisme 
 
 (4.2.1)  \hspace{3cm}  ${\text{D}}_{q}/I \stackrel \simeq  \longrightarrow {\text{End}}_{R[x^{p}]}(R[x]).$  \\
 
  {\it{D\'emonstration}}. La $R[x^{p}]$-lin\'earit\'e r\'esulte du lemme 3.1. Comme les deux modules ont m\^eme rang $p$ sur $R[x^{p}]$, il suffit de prouver l'injectivit\'e. Or, si l'on choisit un repr\'esentant   $P=\sum _{i=0}^{p-1}a_{i}(x)\delta^{i} \in {\text{D}}_{q}$  de $\overline{P} \in {\text{D}}_{q}/I$   et si celui-ci  v\'erifie $P(x^{j})=0$ pour tout $j\in[0,p-1]$, on a clairement $P=0$.

Revenons  \`a la d\'emonstration de la proposition :   les deux membres de (4.1.3)  \'etant  de m\^eme rang $p^{2}$ sur  $\widehat {{\text{Z}}_{q}}$, il suffit de v\'erifier la surjectivit\'e. Comme ils sont de plus complets,   gr\^ace au lemme de Nakayama, on est ramen\'e \`a prouver la surjectivit\'e de la   r\'eduction modulo $I$ de (4.1.3). Cette r\'eduction s'identifie canoniquement \`a l'isomorphisme (4.2.1) : d'o\`u la proposition. $\square$

{\it{Remarque}}.   Comme dans \cite{gqlst}, (p.19 warning) on notera que l'application $\Phi : {\text{Z}}_{q} \rightarrow {\text{Z}}_{q} \simeq {\text{Z}}( {\text{End}}_{ {\text{Z}} _{q}}({\text{Z}}(R[x])))$ (${\text{Z}}$ d\'esignant le centre dans le second membre)  ;  $\delta^{p} \rightarrow D^{p}$,   $x^{p} \rightarrow X^{p}$  n'est pas l'identit\'e  (mais induit  un automorphisme  de $\widehat{ {\text{Z}} _{q}}$). Nous ignorons si l'extension $\Phi : {\text{Z}}_{q} \rightarrow {\text{Z}}_{q}$ permet de neutraliser ${\text{D}}_{q}$ au-dessus de ${\cal{U}}$.

{\bf{Corollaire 4.3}}. La cat\'egorie des $\widehat{ {\text{D}}_{q}}$-modules est  \'equivalente \`a celle des $\widehat{{\text{Z}}_{q}}$-modules.     
  
{\it{D\'emonstration}}. C'est simplement l'\'equivalence de Morita qui \`a un $\widehat{ {\text{D}}_{q}}$-module $E$ associe le $\widehat{{\text{Z}}_{q}}$-module 
${\text{Hom}}_{\widehat{ {\text{D}}_{q}}}(\widehat{ {\text{Z}}(R[x])}, E)$  et \`a un $\widehat{{\text{Z}}_{q}}$-module $F$ associe le $\widehat{ {\text{D}}_{q}}$-module $\widehat{ {\text{Z}}(R[x])}\otimes_{ \widehat{{{\text{Z}}_{q}}}} F$.$\square$

Notons maintenant, pour la suite,  qu'on dispose d'un isomorphisme $\widehat {{\text{Z}}_{q}} \simeq R[x][[\xi]]$ ; $x^p \mapsto x , \delta^p\mapsto \xi$.

{\bf{D\'efinition 4.4}}. Un $R[x]$-\emph{module de Higgs} est un module $H$ muni d'un endomorphisme $R[x]$-lin\'eaire $\xi_{H}$. On dit qu'il est \emph{quasi-nilpotent} si pour tout $h \in H$, il existe $d \geq 0$ tel que $\xi_{H}^d(h) = 0$.

{\bf{D\'efinition 4.5}}. Une \emph{$\sigma$-d\'erivation} sur un $R[x]$-module $M$ est une application $R$-lin\'eaire $D_{M} : M \to M$ telle que $D_{M}(fm) = \delta(f)m + \sigma(f)D_{M}(m)$ pour tous $f\in R[x]$, $m\in M$.
On dit que celle-ci est \emph{quasi-nilpotente} si pour tout $m \in M$, il existe un entier naturel $d$ tel que $D_{M}^d(m)= 0$.

On a alors le r\'esultat suivant, qui ne fait que paraphraser le pr\'ec\'edent lorsqu'on se limite aux objets quasi-nilpotents (voir d\'efinition 5.4 et la remarque suivant la proposition 5.5. dans \cite{gqlst}) :

{\bf{Corollaire 4.6}}. La cat\'egorie des $R[x]$ modules $M$ munis d'une $\sigma$-d\'erivation quasi-nilpotente est \'equivalente \`a celle des  $R[x]$-modules de Higgs $H$ quasi-nilpotents.

Les corollaires 4.3 et 4.6 restent \'evidemment vrais lorsque l'on \'etend les scalaires, par exemple de $R$ \`a  ${\overline{{\Bbb{Q}}}}$ (ou tout autre anneau interm\'ediaire), ou  avec des anneaux de s\'eries formelles en $x$ (avec conditions de convergence \'eventuelle) \`a la place de $R[x]$.

{\it{Remarque}}. L'application    $\Phi : {\text{Z}}_{q} \rightarrow {\text{Z}}_{q}$  ci-dessus est naturellement la restriction  de l'application $\Phi : \widehat{ {\text{D}}_{q}} \rightarrow \widehat{ {\text{D}}_{q}}$ d\'efinie par $\Phi(\delta^{i})=D^{i}(1)$ (e.g. $\Phi(\delta)= \frac{x^{p-1}\delta^{p}}{[p-1]!}$). On peut alors reformuler, dans le corollaire 4.6,  le passage de $M$ \`a $H$  comme la consid\'eration des points fixes de $M$ par  $\Phi$ (cf. \cite{gqlst}, Def. 4.9, prop. 5.7 pour le cas analogue en caract\'eristique $p>0$). 

Signalons pour terminer qu'il n'y a aucune difficult\'e \`a \'etendre, pour $n\geq1$,  ces r\'esultats au cas o\`u $q$ est une racine primitive $p^{n}$-i\`eme de l'unit\'e.


\begin{thebibliography}{99} 

  
\bibitem{backelin} {Backelin, E.} :  {\it{Endomorphisms of quantized Weyl algebras}}. arXiv:1007.2620v1.

\bibitem{bmr} {Bezrukavnikov, R. ; Mirkovi�c, I. ;  Rumynin, D.} : {\it{Localization of modules for a semi-simple Lie algebra in prime characteristic (with an  appendix by R. Bezrukavnikov and S. Riche)}}. Ann. of Math. (2) 167 (2008), 945-991.   MR 2415389. 



 
\bibitem{gqlst} {Gros, M. ;    Le Stum, B. ;  Quiros, A.} :  {\it{A Simpson correspondence in positive characteristic}}. Publ. RIMS Kyoto Univ. 46 (2010) 1-35. MR 2662614.

\bibitem{EGA} {Grothendieck, A. ; Dieudonn\'e, J.A.} : {\it{\'El\'ements 
de G\'eom\'etrie Alg\'ebrique I (seconde \'edition).}}  Springer-Verlag (1971).  

 \bibitem{OV} {Ogus, A. ; Vologodsky, V.} : {\it{Nonabelian Hodge theory in characteristic $p$}}, Publ. Math. Inst. Hautes \'Etudes Sci. 106 (2007), 1-138.  MR 2373230.

 

\bibitem{ta} {Tanisaki, T. } :  {\it{Differential operators on quantized flag manifolds at roots of unity II}}.  arXiv:1101.5848v1.

\bibitem{vdp} {van der Put, M.} :    {\it{Differential equations in characteristic $p$}}, Compos. Math. 97 (1995), 227- 
251.  MR 1355126. 

 

\end{thebibliography}
\end{document}